\documentclass[11pt,reqno]{article}
\usepackage{amsmath,mathrsfs,graphicx,amssymb,amsfonts,color}
\usepackage{cite}

\font\bg=cmbx10 scaled\magstep1

\font\small=cmr8

\newtheorem{newlemma}{{\bf Lemma}}

\newtheorem{newteorem}{{\bf Theorem}}

\newenvironment{teorem}{\begin{newteorem}{\hspace{-0.5
em}{\bf.}}}{\end{newteorem}}

\newtheorem{newkorolari}{{\bf Corollary}}

\newtheorem{newdefine}{{\bf Definition}}

\newenvironment{define}{\begin{newdefine}{\hspace{-0.5
em}{\bf.}}}{\end{newdefine}}

\newtheorem{newquestion}{{\bf Question}}

\newtheorem{newkonjek}{{\bf Conjecture}}

\newtheorem{newexample}{{\bf Example}}

\begin{document}
\tolerance=10000
\baselineskip18truept
\newbox\thebox
\global\setbox\thebox=\vbox to 0.2truecm{\hsize 0.15truecm\noindent\hfill}
\def\boxit#1{\vbox{\hrule\hbox{\vrule\kern0pt
\vbox{\kern0pt#1\kern0pt}\kern0pt\vrule}\hrule}}
\def\qed{\lower0.1cm\hbox{\noindent \boxit{\copy\thebox}}\bigskip}
\def\ss{\smallskip}
\def\ms{\medskip}
\def\bs{\bigskip}
\def\c{\centerline}
\def\nt{\noindent}
\def\ul{\underline}
\def\ol{\overline}
\def\lc{\lceil}
\def\rc{\rceil}
\def\lf{\lfloor}
\def\rf{\rfloor}
\def\ov{\over}
\def\t{\tau}
\def\th{\theta}
\def\k{\kappa}
\def\l{\lambda}
\def\L{\Lambda}
\def\g{\gamma}
\def\d{\delta}
\def\D{\Delta}
\def\e{\epsilon}
\def\lg{\langle}
\def\rg{\rangle}
\def\p{\prime}
\def\sg{\sigma}
\def\ch{\choose}

\newcommand{\ben}{\begin{enumerate}}
\newcommand{\een}{\end{enumerate}}
\newcommand{\bit}{\begin{itemize}}
\newcommand{\eit}{\end{itemize}}
\newcommand{\bea}{\begin{eqnarray*}}
\newcommand{\eea}{\end{eqnarray*}}
\newcommand{\bear}{\begin{eqnarray}}
\newcommand{\eear}{\end{eqnarray}}

\centerline{\Large \bf Some families of graphs with no nonzero real domination roots}
\bigskip

\bs

\baselineskip12truept
\centerline{S. Jahari$^{}${}\footnote{\baselineskip12truept\it\small 
Corresponding author. E-mail: s.jahari@gmail.com} and S. Alikhani }
\baselineskip20truept
\centerline{\it Department of Mathematics, Yazd University}
\vskip-8truept
\centerline{\it  89195-741, Yazd, Iran}

\vskip-0.2truecm
\nt\rule{16cm}{0.1mm}

\nt{\bg ABSTRACT}
\medskip

\baselineskip14truept

\nt{Let $G$ be a simple graph of order $n$.
The {\em domination polynomial} of $G$ is the polynomial
$D(G, x)=\sum_{i=\gamma(G)}^{n} d(G,i) x^{i}$, where $d(G,i)$ is the number of dominating sets of $G$ of size $i$ and $\gamma(G)$ is the domination number of $G$. A root of $D(G, x)$ is called a domination root of $G$. Obviously, $0$ is a domination root of every graph $G$ with multiplicity $\gamma(G)$. In the study of the domination roots of graphs, this naturally raises
the question: Which graphs have no nonzero real domination roots?  In this paper we present some families of graphs whose have this property.  }

\ms

\nt{\bf Mathematics Subject Classification:} {\small 05C31, 05C60.}

\nt\rule{16cm}{0.1mm}

\baselineskip20truept

\section{Introduction}

\nt All graphs in this paper are simple of finite orders, i.e., graphs are undirected with no loops or
parallel edges and with finite number of vertices.  Let $G=(V,E)$ be a simple graph.
For any vertex $v\in V(G)$, the {\it open neighborhood} of $v$ is the
set $N(v)=\{u \in V (G) | uv\in E(G)\}$ and the {\it closed neighborhood} of $v$
is the set $N[v]=N(v)\cup \{v\}$. For a set $S\subseteq V(G)$, the open
neighborhood of $S$ is $N(S)=\bigcup_{v\in S} N(v)$ and the closed neighborhood of $S$
is $N[S]=N(S)\cup S$.

\nt The {\it complement} $G^c$ of a graph $G$ is a graph with the same vertex set as $G$ and with the property that two vertices are adjacent in $G^c$ if and only if they are not adjacent in $G$.

\nt A set $S\subseteq V(G)$ is a {\it dominating set} if $N[S]=V$ or equivalently,
every vertex in $V(G)\backslash S$ is adjacent to at least one vertex in $S$.
The {\it domination number} $\gamma(G)$ is the minimum cardinality of a dominating set in $G$.
For a detailed treatment of domination theory, the reader is referred to~\cite{domination}.

\nt Let ${\cal D}(G,i)$ be the family of dominating sets of a graph $G$ with cardinality $i$ and
let $d(G,i)=|{\cal D}(G,i)|$.
The {\it domination polynomial} $D(G,x)$ of $G$ is defined as
${D(G,x)=\sum_{ i=\gamma(G)}^{|V(G)|} d(G,i) x^{i}}$ (see \cite{euro,saeid1,kotek}); this polynomial is the generating polynomial for the number of dominating sets of each cardinality.
Similar to generating polynomials for other combinatorial sequences, such as independents sets in a graph 
\cite{bn,brownhickmannowa,chudnovsky,fish90,gutman,gut90,gut92}, have
 attracted recent attention, to name but a few. The algebraic encoding of salient counting sequences allows one to not only develop formulas more easily, but also, often, to prove unimodality results via the nature of the the roots of the associated polynomials (a well known result of Newton states that if a real polynomial with positive coefficients has all real roots, then the coefficients form a unimodal sequence (see, for example, \cite{comtet}). A root of $D(G, x)$ is called a {\it domination root} of $G$ (see \cite{few,brown}). The set of distinct non-zero roots  of $D(G, x)$ is denoted by $Z^*(D(G, x))$. It is known that $-1$ is \underline{not} a domination root as the number of dominating sets in a graph is always odd \cite{brouwer}. On the other hand, of course, $0$ is a domination root of every graph $G$ with multiplicity $\gamma(G)$. The existing research on the roots of domination
polynomials has been restricted to those graphs with exactly two, three or exactly
four domination roots (\cite{euro,four}). Also in \cite{brown} Brown and Tufts
studied the location of the roots of domination polynomials for some families of graphs such as
bipartite cocktail party graphs and complete bipartite graphs. In particular, they showed that the
set of all domination roots is dense in the complex plane.

 \nt In the study of the domination roots of graphs, this naturally raises
the question: Which graphs have no nonzero real domination roots? In this paper we would like to present some families of graphs with this property.
Let ${\cal G}$ be the family of graphs and $\mathcal{CG}=\big\{G\in {\cal G}|Z^*(D(G,x))\subseteq  \mathbb{C}\big\}$.

 \nt In the next section we present some families of graphs whose are in $\mathcal{CG}$.  In Section $3$ we consider  the complement of the friendship graphs, $F_n^c$ and compute their domination polynomials, exploring the nature and location of their roots. As a consequence we show that $F_n^c\in \mathcal{CG}$.

\section{Some families of graphs in $\cal{CG}$}

\nt In the beginning of the  study of domination roots of graphs, one can see that there are graphs with no nonzero real domination roots except zero. As examples, the complete graph $K_n$ for odd $n$ and the complete bipartite graph $K_{n,n}$ for even $n$, are in $\cal{CG}$. With these motivation, in \rm\cite{SSM,few} the authors asked the question: ``Which graphs have no nonzero real domination roots?" In other words,
which graph lie in $\cal{CG}$?.

\nt In this section we use the existing results on domination polynomials to find some families of graphs whose are in $\cal{CG}$.  We need some preliminaries.

\nt The join $G = G_1+ G_2$ of two graphs $G_1$ and $G_2$ with disjoint vertex sets $V_1$ and $V_2$ and edge sets $E_1$ and $E_2$ is the graph union $G_1\cup G_2$ together with all the edges joining $V_1$ and $V_2$.
\nt The following theorem gives a formula for the domination polynomial of join of two graphs.

\begin{teorem}\label{theorem8}{\rm \cite{euro}}
Let $G$ and $H$ be nonempty graphs of order $n$ and $m$, respectively. Then,
\begin{eqnarray*}
D(G + H,x) = ((1 +x)^n -1) ((1 + x)^m -1) + D(G, x) + D(H, x).
\end{eqnarray*}
\end{teorem}

\nt For two graphs $G = (V,E)$ and $H=(W,F)$, the corona $G\circ H$ is the graph arising from the
disjoint union of $G$ with $| V |$ copies of $H$, by adding edges between
the $i$th vertex of $G$ and all vertices of $i$th copy of $H$ \cite{Fruc}.  We need the following theorem which is for computing  the domination
polynomial of the corona products of two graphs.

\begin{teorem}\label{theorem7}{\rm \cite{Oper,Kot}}
Let $G = (V,E)$ and $H=(W,F)$ be nonempty graphs of order $n$ and $m$, respectively. Then
\begin{eqnarray*}
D(G\circ H,x) = (x(1 + x)^m + D(H, x))^n.
\end{eqnarray*}
\end{teorem}


\nt To present some families of graphs in $\mathcal{CG}$, we recall the existing results.

\nt A $k$-star, $S_{k,n-k}$, has vertex set $\{v_1, \ldots , v_n\}$ where $<\{v_1, v_2, \ldots , v_k\}>  \cong K_k$ and $N(v_i) = \{v_1, \ldots, v_k\}$ for $k + 1 \leq i \leq n$.

\nt The book graph $B_n$ can be constructed by bonding $n$ copies of the cycle graph $C_4$ along a common edge $\{u, v\}$.  In \rm\cite{jason} it was proved that, for every $n \in \mathbb{N}$,
\begin{eqnarray*}
 D(B_n,x)=(x^2+2 x)^n(2x+1) + x^2(x+1)^{2n} - 2x^n.
 \end{eqnarray*}
The following theorem gives  some families of graphs whose are in $\cal{CG}$.

\begin{teorem}
\begin{enumerate}
\item[(i)] {\rm\cite{jason}} Every graph $H$ in
the family $\{G\circ K_{2n}, (G\circ K_{2n})\circ K_{2n}, ((G\circ K_{2n})\circ K_{2n})\circ K_{2n},\cdots \}$ lie in $\cal{CG}$.

\item[(ii)] {\rm \cite{jahari}}    For odd   $n$ and even $k$, the $k$-star $S_{k,n-k}$ is in $\cal{CG}$.
\item[(iii)] {\rm \cite{jahari}}  For odd   $n$ and odd $k$, every graph $H$ in the family
$$\{G\circ S_{k, n-k}, (G\circ S_{k, n-k})\circ S_{k, n-k}, ((G\circ S_{k, n-k})\circ S_{k, n-k})\circ S_{k, n-k},\cdots \}$$ lie in  $ \cal{CG}$.
\item[(iv)] \rm\cite{jason} Every graph $H$ in
the family $\{G\circ B_2, (G\circ B_2)\circ B_2, ((G\circ B_2)\circ B_2)\circ B_2,\cdots \}$ lie in  $ \cal{CG}$.

\end{enumerate}
\end{teorem}

\begin{figure}[!ht]
\hspace{1.9cm}
\includegraphics[width=11cm,height=4.1cm]{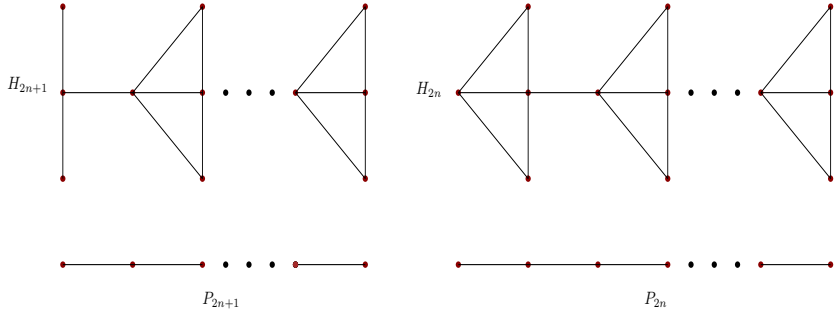}
\caption{ \label{figure1} Graphs $H_{2n+1}$ and $H_{2n}$, respectively. }
\end{figure}

\nt In \cite{Lev},  Levit  and Mandrescu constructed a family of graphs $H_n$ from the path $P_n$ by
the ``clique cover construction'', as shown in Figure \ref{figure1}. By $H_0$ we mean the null graph.

\nt The following theorem  gives formula for the domination polynomials of $H_n$ graphs:

\begin{teorem}\label{theorem9}\rm\cite{aj}
 Let   $H_n$ be the  graphs in the Figure \ref{figure1}.
 \begin{enumerate}
 \item[(i)] For every $n \in \mathbb{N}$, $D(H_{2n},x)=(x^4+4x^3+6x^2+2x)^n.$
 \item[(ii)]  For every $n \in \mathbb{N}$, $D(H_{2n+1},x)=(x^3+3x^2+x)(x^4+4x^3+6x^2+2x)^n.$
 \end{enumerate}
\end{teorem}

\nt Here using Theorem \ref{theorem9} we present another families of graphs in $\cal{CG}$.

\begin{teorem}
\begin{enumerate}
\item[(i)] The graphs of the form $H_n+ H_n$,  $H_{n+1}+ B_n$, for $n\geq 3$, and the graphs of the form $B_n+B_n$, for  odd  $n$ are in $\cal{CG}$.
\item[(ii)]
 The graphs of the form $B_{n+1}+ B_n$, for  even $n$, and $B_{n+1}+ H_n$, for $n\geq 4$ are in $\cal{CG}$.
\end{enumerate}
\end{teorem}

\nt{\bf Proof.}
Since  the coefficients of domination polynomials are positive integers,  we investigate domination roots for $x\leq 0$.
\begin{enumerate}
\item[(i)]
 By theorem~\ref{theorem8} we can deduce that for each natural number $n\geq 3$,
\[
D(H_n+ H_n,x)= ((1+x)^{|V(H_n)|}-1)^2+2D(H_{n},x).
\]
To obtain the domination roots of $H_n+ H_n$, we shall solve the following equation:
 \begin{eqnarray}\label{eq}
  ((1+x)^{|V(H_n)|}-1)^{2} = - 2D(H_{n},x).
 \end{eqnarray}
 We consider two cases, and show that in each there is no nonzero solution.
\begin{itemize}
\item If $n\geq 3$ is even, i.e., $n=2k$ for some $k\in \mathbb{N}$.  Then the equation (\ref{eq}) is equivalent to the following equation
 \begin{eqnarray}
  ((1+x)^{4k}-1)^{2} &=& - 2(x^4+4x^3+6x^2+2x)^k\nonumber\\
  &=&-2((1+x)^4-2x-1)^k.
 \end{eqnarray}
For $x\leq 0$,  the above equality is true
just for real number $0$. Because for nonzero real number the left
side of this equality is positive but the right side is negative.

\item If $n\geq 3$ is odd, i.e., $n=2k+1$, $n=2k$ for some $k\in \mathbb{N}$. Then the equation (\ref{eq}) is is equivalent to the following equation
 \begin{eqnarray}\label{eq2}
  ((1+x)^{4k+3}-1)^{2} &=& - 2(x^3+3x^2+x)(x^4+4x^3+6x^2+2x)^k\nonumber\\
  &=&-2((1+x)^3-2x-1)((1+x)^4-2x-1)^k.
 \end{eqnarray}
We consider the following different cases, and show in each there is no nonzero solution.
 If  $x \leq -1$,  there are no real solutions $x$.  Because, it is easy to see that for $-2\leq x\leq -1$, the left side of \ref{eq2} is positive but its  right side
 is negative. Also for $x<-2$,   the left side of equality (\ref{eq2}) is greater than the right side. Now suppose that $-1 < x < 0$.
\begin{enumerate}
\item If $k$ is even and $-\frac{1}{2}\leq x < 0,$   the left side of equality (\ref{eq2}) is greater than the right side, a contradiction.
\item If $k$ is odd and $-\frac{1}{2}\leq x < 0,$   the left side of equality (\ref{eq2}) is positive but the right side is negative, a contradiction.
\item For every $k$ and $-1 < x <- \frac{1}{2},$  there are no real solutions $x$. Because  the left side of equality (\ref{eq2}) is positive but the right side is negative.
\end{enumerate}
 \end{itemize}
 \nt The other cases are similar to this case.

\item[(ii)] It is similar to proof of Part (i).\quad\qed
\end{enumerate}

\begin{figure}[ht]
\hspace{5.cm}
\includegraphics[width=5.5cm,height=4.2cm]{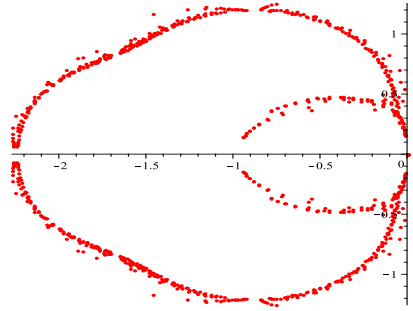}
\caption{\label{figure5} Domination roots of~$H_n+H_n$, ~for~odd $n$ and $3\leq n\leq 20$.}
\end{figure}

\nt Domination roots of the graphs $H_n+H_n$, for odd $n$ and~ $3\leq n\leq 20$ has shown in Figure \ref{figure5}.

\section{Domination roots of the complement of the friendship graphs}

\begin{figure}
\begin{center}
\includegraphics[width=6in]{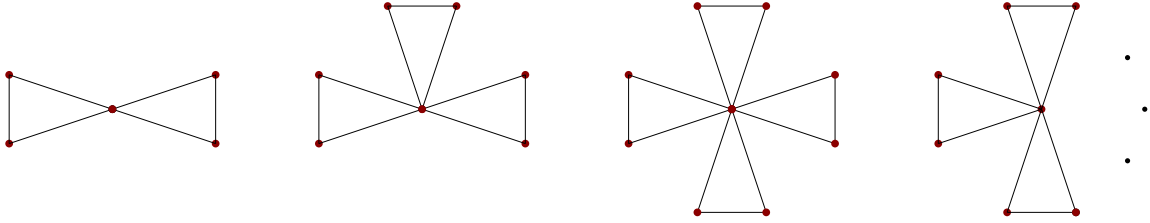}
\caption{Friendship graphs $F_2, F_3, F_4$ and $F_n$, respectively.}
\label{Dutch}
\end{center}
\end{figure}

\nt The friendship (or Dutch-Windmill) graph $F_n$ is a graph that can be constructed by coalescence $n$
copies of the cycle graph $C_3$ of length $3$ with a common vertex. The Friendship Theorem of Paul Erd\"{o}s,
Alfred R\'{e}nyi and Vera T. S\'{o}s \cite{erdos}, states that graphs with the property that every two vertices have
exactly one neighbour in common are exactly the friendship graphs.
Figure \ref{Dutch} shows some examples of friendship graphs.

\nt The following theorem states that for each odd $n$, the friendship graph  $F_{n}$ lie in $\cal{CG}$.
\begin{teorem}{\rm \cite{jason}}
\begin{enumerate}
\item[(i)]
  For every $n\in \mathbb{N}$, $D(F_{n},x) = (2x + x^2)^n +x(1 + x)^{2n}.$
\item[(ii)] For odd $n$, $F_{n}\in \cal{CG}$.
\end{enumerate}
\end{teorem}

 \nt Domination polynomials, exploring the nature and location of domination  roots of friendship graphs has studied in \cite{jason}.  It is natural to ask about the domination polynomial and the domination roots of the complement of the friendship graphs.

\nt The Tur\'{a}n graph $T(n,r)$ is a complete multipartite graph formed by partitioning a set of $n$ vertices into $r$ subsets, with sizes as equal as possible, and connecting two vertices by an edge whenever they belong to different subsets. The graph will have $(n~ mod ~ r)$ subsets of size $\lceil\frac{n}{r}\rceil$, and $r - (n~ mod ~ r)$ subsets of size $\lfloor\frac{n}{r}\rfloor$. That is, it is a complete $r$-partite graph
\[ K_{\lceil\frac{n}{r}\rceil, \lceil\frac{n}{r}\rceil,\ldots,\lfloor\frac{n}{r}\rfloor,\lfloor\frac{n}{r}\rfloor}.\]
\nt The Tur\'{a}n graph $T(2n,n)$ can be formed by removing  a perfect matching,  $n$ edges no two of which are adjacent,  from a complete graph $K_{2n}$. As Roberts (1969) showed, this graph has boxicity exactly $n$; it is sometimes known as the Roberts graph \cite{Robert}. If $n$ couples go to a party, and each person shakes hands with every person except his or her partner, then this graph describes the set of handshakes that take place; for this reason it is also called the cocktail party graph. So, the cocktail party graph $CP(t)$ of order $2t$ is the graph with vertices $b_1, b_2, \cdots, b_{2t}$
in which each pair of distinct vertices form an edge with the exception of the pairs
$\{b_1 , b_2 \}, \{b_3 , b_4\}, \ldots, \{b_{2t- 1}, b_{2t}\}$. 

\begin{figure}
\begin{center}
\includegraphics[width=6in]{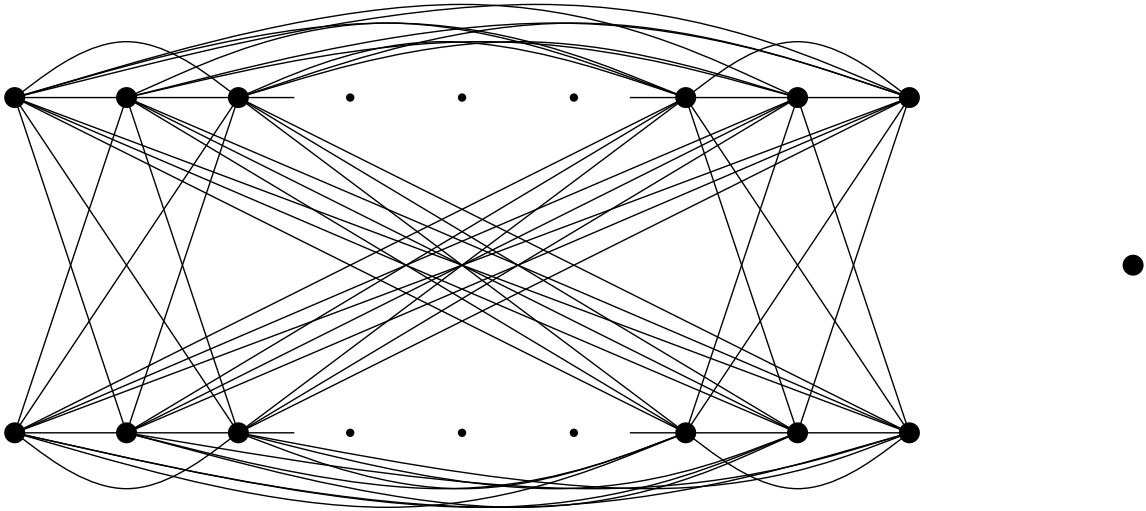}
\caption{Complement of the friendship graph $F_n$.}
\label{CDutch}
\end{center}
\end{figure}

\nt It is easy to check that the complement of the  friendship graph $F_n$ is $CP(n)\cup K_1$.
Figure \ref{CDutch} shows  the complement of the  friendship graph $F_n$.

\begin{teorem}\label{d.p.dutch}
For every $n\in \mathbb{N}$, $D(F_n^c,x) = \big((1+ x)^{2n} -(1 + 2nx)\big)x.$
\end{teorem}
\nt{\bf Proof.}  An elementary observation is that if $G_1$ and $G_2$ are graphs of orders $n_1$ and $n_2$,
respectively, then
\[ D(G_1 \cup G_2,x) = D(G_1, x) D(G_2, x).\]
Clearly $D(K_{1},x) = x$ and there are no dominating sets of size $1$ in $CP(n)$. Therefore $$D(CP(n),x) = (1+ x)^{2n} -(1 + 2nx).\quad\qed$$

\begin{figure}[ht]
\hspace{4cm}
\includegraphics[width=7cm]{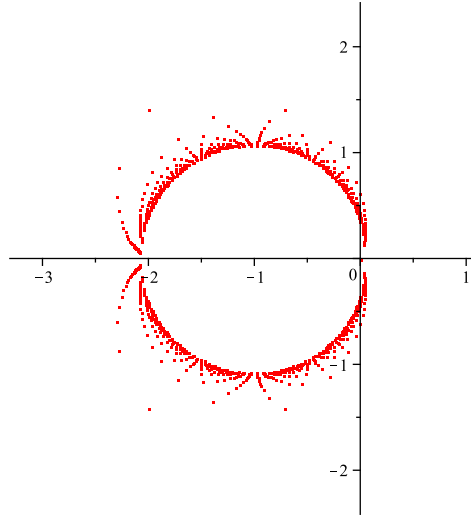}
\caption{\label{figure2'} Domination roots of graphs $F_n^c$, for $1 \leq n \leq 30$.}
\end{figure}

\nt In \cite{brown} a family of graphs was produced with roots just barely in the right-half plane (showing that not all domination polynomials are stable), but Figure~\ref{figure2'}  provides an explicit family (namely the $F_n^c$) whose domination roots have positive real part.
\ms

\nt The domination roots of complement of the  friendship graphs exhibit a number of interesting properties (see Figure~\ref{figure2'}). Even though we cannot find the roots explicitly, there is much we can say about them.


\nt Here we  prove that  for each natural number $n$,  the complement of  the friendship graphs $F_{n}^c$ lie in $\cal{CG}$.

\begin{teorem}
For every natural number $n$,  the complement of  the friendship graphs $F_{n}^c$ lie in $\cal{CG}$.
     \end{teorem}
\nt{\bf Proof.}  It's suffices to show that for each natural $n$,  the cocktail party graph $CP(n)$ is in $\cal{CG}$. By Theorem \ref{d.p.dutch}, for every $n\in
\mathbb{N}$, $D(CP(n),x) = (1+ x)^{2n} -(1 + 2nx). $ If
$D(CP(n),x)=0$  then for $x \neq 0$, we have
\[
(1+ x)^{2n} =1 + 2nx.
\]
We consider three cases, and show in each there is no nonzero solution.
\begin{itemize}
\item  $x > 0:$  Obviously the above equality is true
just for real number $0$, since for nonzero real number the left
side of equality is greater than the right side.
\item $x \leq -1:$ In this case the left side is greater than $0$
and the right side $1 + 2nx$ is less
than $-1$, a contradiction.
\item  $-1 < x < 0:$  In this case obviously there are no real solutions $x$,  the left side of equality is greater than   the right side.
\end{itemize}

\nt Thus in any event, there are no nonzero real domination roots of the cocktail party graph.
\hfill \qed




\nt The plot in Figure~\ref{figure2'} suggests that the roots tend to lie on a curve. In order to find the limiting curve, we will need a definition and a well known result.

\begin{define}
If ${f_n(x)}$ is a family of (complex) polynomials, we say that a number $z \in \mathbb{C}$ is a limit of roots of ${f_n(x)}$ if either $f_n(z) = 0$ for all sufficiently large $n$ or z is a limit point of the set $\mathbb{R}({f_n(x)})$, where $\mathbb{R}({f_n(x)})$ is the union of the roots of the $f_n(x)$.
\end{define}

\nt The following restatement of the Beraha-Kahane-Weiss theorem \cite{bkw}  can be found in \cite{brownhickman}.

\begin{teorem}\label{bkw}
Suppose ${f_n(x)}$ is a family of polynomials such that
\begin{eqnarray}
f_n(x) = \alpha_1(x)\lambda_1(x)^n + \alpha_2(x)\lambda_2(x)^n + ... + \alpha_k(x)\lambda_k(x)^n
\end{eqnarray}
where the $\alpha_i(x)$ and the $\lambda_i(x)$ are fixed non-zero polynomials, such that for no pair $i \neq j$ is $\lambda_i(x) \equiv \omega\lambda_j(x)$ for some $\omega \in \mathbb{C}$ of unit modulus. Then $z \in \mathbb{C}$ is a limit of roots of ${f_n(x)}$ if and only if either
\begin{itemize}
\item[(i)] two or more of the $\lambda_i(z)$ are of equal modulus, and strictly greater (in modulus) than the others; or
\item[(ii)] for some $j$, $\lambda_j(z)$ has modulus strictly greater than all the other $\lambda_i(z)$, and $\alpha_j(z) = 0$
\end{itemize}
\end{teorem}

\nt The following Theorem gives the limits of the domination roots of  $F_n^c$.

\begin{teorem}
The limit of domination roots of $F_n^c$  is the unit circle with center $-1$.
\end{teorem}

\nt{\bf Proof.} By Theorem~\ref{d.p.dutch}, the domination polynomial of $F_n^c$ is,
\begin{eqnarray*}
D(F_n^c,x) &=& x((x+1)^2)^{n} - x(1 + 2nx)\\
&=&\alpha_1(x)\lambda_1^n(x) + \alpha_2(x)\lambda_2^n(x),
\end{eqnarray*}
\nt where
\[ \alpha_{1}(x) = x,~~\lambda_{1}(x) = (x+1)^{2},\]
and
\[\alpha_{2}(x) =x+2nx^2,~~\lambda_{2}(x) = 1.\]

\nt Clearly there is no $\omega \in \mathbb{C}$ of modulus $1$ for which $\lambda_{1} = \omega \lambda_{2}$ (or vice versa). Also, $\alpha_1, \mbox{ and } \alpha_2$ are not identically zero. Therefore, the initial conditions of Theorem~\ref{bkw} are satisfied. Now,  $|x-(-1)|^2=1 $ implies that $x$ lies on the circle centred at $-1$.\quad\qed

\nt{\bf Conclusion.} In this paper we presented some families of graphs whose non-zero domination roots are complex. We think that 
these kind of graphs shall have specific geometrical properties.   However,
until now all attempts to find these properties failed, and it remains as open problem. 


\end{document}